\documentclass[11pt]{article}
\usepackage{graphicx}
\usepackage{amsthm, amsmath, amssymb, tikz, bm}
\usepackage{mathtools, intcalc}
\usepackage{xifthen}
\usepackage{comment}
\usepackage[colorlinks=true,
linkcolor=blue,citecolor=blue,
urlcolor=blue]{hyperref}

\usepackage[margin=1.2in]{geometry} 

\usepackage[shortlabels]{enumitem}
\usepackage{todonotes}
\usetikzlibrary{calc,shapes, backgrounds}
\allowdisplaybreaks

\usepackage[utf8]{inputenc}
\usepackage[T1,T2A]{fontenc}
\usepackage[russian,english]{babel}

\newtheorem{lemma}{Lemma}
\newtheorem{theorem}{Theorem}

\newtheorem{proposition}{Proposition}

\usepackage[normalem]{ulem}

\title{On total bondage number of graphs}

\author{
  E.G.K.M.Gamlath\thanks{University of Mississippi, University, MS 38677, (\texttt{egkmgamlath@gmail.com})}
  \and
  Bing Wei\footnotemark[1]
}

\begin{document}
\maketitle
 
\begin{abstract}
In this paper, we explore the concept of total bondage in finite graphs without isolated vertices. A vertex set $D$ is considered a total dominating set if every vertex $v$ in the graph $G$ has a neighbor in $D$. The minimum cardinality of all total dominating sets in $G$ is denoted as $\gamma_t(G)$. A total bondage edge set $B$ is a subset of the edges of $G$ such that the removal of $B$ from $G$ does not create isolated vertices, and the total dominating number of the resulting graph $G-B$ is strictly greater than $\gamma_t(G)$. The total bondage number of $G$, denoted $b_t(G)$, is defined as the minimum cardinality of such total bondage edge sets. Our paper establishes upper bounds on $b_t(G)$ based on the maximum degree of a graph. Notably, for planar graphs with minimum degree $\delta(G) \geq 3$, we prove $b_t(G) \leq \Delta + 8$ or $b_t(G) \leq 10$. Additionally, for a connected planar graph with $\delta(G) \geq 3$ and $g(G) \geq 4$, we show that $b_t(G) \leq \Delta + 3$ if $G$ does not contain an edge with degree sum at most 7. We also improve some upper bounds of the total bondage number for trees, enhance existing lemmas, and find upper bounds for total bondage in specific graph classes.
\end{abstract}

\section{Introduction}
Graph theory serves as an invaluable tool for modeling and understanding complex relationships within various domains. In this paper, we delve into the concept of the total bondage number, an intriguing parameter that sheds light on the structural vulnerability of graphs. The concept of the total bondage number was introduced by Kulli and Patwari \cite{kulli1991advances} in 1991, and denoted as \(b_t(G)\), provides a quantitative measure of the minimal edge removal required to increase the total domination of a graph \(G\) without isolated vertices.

The total domination number itself is a well-studied parameter. For a finite graph $G$ without isolated vertices, a vertex set $D$ of $G$ is said to be a total dominating set of $G$, if every vertex $v\in V(G)$ has a neighbor in $D$. Define $\gamma_t (G)$ to be the minimum cardinality among all total dominating sets of $G$. Hu and Xu demonstrated in \cite{hu2011complexity} showed that determining the total bondage number for general graphs is an NP-hard problem. 

In this paper we prove that,
\begin{theorem}
    \label{totalbondagemindegree3planar}
For a connected planar graph with $\delta(G)\geq 3$,
    \begin{center}
    
$b_t(G)\leq \Delta(G) + 8
$ or $b_t(G)\leq 10.  $
\end{center}

\end{theorem}
Using some lemmas and some configurations of a planar graph with $\delta(G)\geq 3$ (Theorem \ref{borodin}), we prove the theorem \ref{totalbondagemindegree3planar}.

In \cite{rad2014some}, it was demonstrated that for any tree $T$ with $\Delta(G)\geq 3$, the total bondage number satisfies $b_t(G)\leq \Delta -1$. Additionally, for a connected graph, the upper bounds were established as $b_t(G)\leq 3\Delta -4$ or $b_t(G) \leq 3\Delta -5$. In this paper, we build upon these results, providing improvements to the aforementioned upper bounds. Furthermore, we enhance the upper bound previously proposed by Sridharan in \cite{sridharan2007total}. Specifically, we establish that, \\

\begin{theorem}
\label{sridharanimprovement}
For a connected tree $T$, excluding cases where $T$ is isomorphic to $P_4$ (path with 4 vertices), $K_{1,3}$ (star with one central vertex connected to three peripheral vertices), or the tree $T^1$ with 7 vertices obtained by subdividing a single edge of $K_{1,3}$ three times, and given that the minimum degree ($\Delta$) of $T$ is greater than or equal to 3, we have $b_t(T) \leq \frac{n-2}{3}$, where $n$ is the number of vertices in $T$.

\end{theorem}


We also establish the following results:

\begin{theorem}
\label{girth4totalbondage}
For a connected planar graph $G$ with $\delta(G) \geq 3$ and girth $g(G) \geq 4$, without an edge with degree sum at most 7, we have $b_t(G) \leq \Delta + 3$.
\end{theorem}

To prove Theorem \ref{girth4totalbondage}, we employ the following lemma:

\begin{lemma}
\label{3vertices at distance 2}
For a graph with $\delta(G) \geq 3$, if there exists a path with two 3-vertices at a distance exactly 2, then $b_t(G) \leq \Delta + 3$.
\end{lemma}

Then, we prove the following theorem.

\begin{theorem}
\label{distant2totalbondage}
For a connected graph $G$ with $\delta(G) \geq 2$, if there exists a path with two 2-vertices at a distance of at most 3, then $b_t(G) \leq \Delta + 1$.
\end{theorem}



\vspace{0.3cm}
To prove the theorem \ref{girth4totalbondage}, we prove the following theorem using the discharging method,

\begin{theorem}
\label{configurationgirth4}
Let $G$ be a connected planar graph with $\delta(G) \geq 3$ and girth $g(G)\geq 4$. Then, $G$ has at least one of the following configurations:
\begin{enumerate}
   \item[(a)] $(3,4^-)$-edge.
   \item[(b)] $v \in S_5$ such that $|N(v)\cap S_3|\geq 4$.
\end{enumerate}
\end{theorem}
Then, we prove the upper bound for total bondage number for complete multipartite graphs.
\begin{theorem}
\label{completemultipartitegraphtotalbondate}
 Let $K_{n_1,n_2,\dots, n_k}$ be a complete multipartite graph with $2\leq n_k\leq n_{k-1}\leq \dots \leq n_2\leq n_1 $. Then, 
 $$b_t(K_{n_1,n_2,\dots, n_k})\leq 4n-2n_1-2$$, where $n$ is the number of vertices in the complete multipartite graph.
\end{theorem}

This paper is organized as follows. In section 2, we give definitions and known results. In section 3, we prove the lemmas and theorems, and in section 4 we provide the graphs with total bondage number 3,4, and 5.

\section{Definitions and known results}

Let $G$ be a connected finite graph with sets of vertices $V(G)$ and edges $E(G)$. For a vertex $v \in V(G)$, $N_G(v)$ represents the open neighborhood of $v$. When the graph $G$ is explicitly known, we simply use $N(v)$ for the open neighborhood of $v$. The degree of a vertex $v$ is denoted by $d(v) = |N(v)|$. The maximum and minimum degrees of $G$ are represented by $\Delta(G)$ and $\delta(G)$, respectively.

Let $S_1$ be the set of pendant vertices of $G$. A vertex $v$ is considered a \textbf{support vertex} if $|N(v)\cap S_1| \neq \emptyset$. $S_n$ represents the set of vertices in $G$ where $d(v) = n$. An $(a, b)$-edge is an edge where the end vertices have degrees $a$ and $b$, respectively. An $(a,b^-)$-edge is an edge where the end vertices have degree $a$ and at most $b$. The distance between vertices $u$ and $v$ is denoted by $d(u,v)$. $C_n$ represents a cycle with $n$ vertices, and an $n$-face is a face with $n$ vertices.A vertex is called an $a^-$-vertex if its degree is at most $a$, and it is referred to as an $a^+$-vertex if its degree is at least $a$. The number of vertices in $G$ is denoted by $|V(G)| = n$. 

The complete bipartite graph $K_{m,n}$ is characterized by two disjoint sets of vertices $V_1$ and $V_2$ with cardinalities $|V_1| = m$ and $|V_2| = n$. Every vertex in $V_1$ is adjacent to every vertex in $V_2$, establishing a comprehensive connectivity pattern.

In contrast, the complete multipartite graph $K_{n_1, n_2, \ldots, n_k}$ consists of $k$ disjoint sets of vertices $V_1, V_2, \ldots, V_k$, each with cardinalities $|V_1| = n_1, |V_2| = n_2, \ldots, |V_k| = n_k$. The key feature here is that each vertex in a set $V_i$ forms an edge with every vertex in all other sets $V_j$ where $j \neq i$.

For a finite graph $G$ without isolated vertices, a vertex set $D$ is a \textbf{total dominating set} of $G$ if every vertex $v\in V(G)$ has a neighbor in $D$. The \textbf{total domination number} is denoted by $\gamma_t (G)$, defined as the minimum cardinality among all total dominating sets of $G$. For any graph $G$, let $TD(G)$ be a minimum cardinality total dominating set of $G$.

A \textbf{total bondage edge set} $B$ is a subset of $E(G)$ such that $G-B$ does not have isolated vertices, and $\gamma_t(G-B)>\gamma_t(G)$. The \textbf{total bondage number} of $G$ is denoted by $b_t(G)$, and $b_t(G) = \min\{|B|: B \text{ is a total bondage edge set of } G\}$ if $G$ has at least one total bondage edge set; otherwise, $b_t(G)=\infty$.


In their work \cite{kulli1991advances}, Kulli and Patwari introduced the total bondage number for a graph. Hu and Xu demonstrated in \cite{hu2011complexity} that determining the total bondage number for general graphs is an NP-hard problem.  Kulli and Patwari \cite{kulli1991advances} computed the exact values of $b_t(G)$ for certain standard graphs, such as a cycle $C_n$ and a path $P_n$ for $n\geq 4$ a complete bipartite graph $K_{m,n}$, and a complete graph $K_n$. Furthermore, Sridharan et al. \cite{sridharan2007total} established that for any positive integer $k$, there exists a tree $T$ with  $b_t(T)=k$.

\begin{proposition} (Kulli et al. \cite{kulli1991advances}) For $n\geq 2$,

\begin{center}
    
$b_t(P_n) = \begin{cases}
        \infty, & \text{if }  n \leq 3\\
        1, & \text{if }  n \geq 4,  \not\equiv 2\ (\textrm{mod}\ 4\\
        2, & \text{if }   n \geq 4,  \equiv 2\ (\textrm{mod}\ 4\\
        \end{cases}$
\end{center}.
    
\end{proposition}

\begin{proposition}(Kulli et al. \cite{kulli1991advances}) For $n\geq 3$,
\begin{center}
    
$b_t(C_n) = \begin{cases}
        \infty, & \text{if }  n = 3\\
        2, & \text{if }  n \geq 4,  \not\equiv 2\ (\textrm{mod}\ 4\\
        3, & \text{if }   n \geq 4,  \equiv 2\ (\textrm{mod}\ 4\\
        \end{cases}$.
\end{center}    
\end{proposition}
In the same paper Kulli et al. \cite{kulli1991advances} prove that, 
\begin{theorem}
\label{kullibipartite}
    Let $K_{m,n}$ be a complete bipartite graph with $2\leq m\leq n$. Then, $$b_t(K_{m,n})=m.$$
\end{theorem}

In \cite{sridharan2007total}, Sridharan et al. established the following upper bounds for the total bondage number of graphs.
\begin{theorem} \cite{sridharan2007total} Let $G$ be a connected graph of order $n\geq 4$. Then,
\begin{enumerate}
    \item $b_t(G)\leq n-1$ if $g(G)\geq 5$,
    \item $b_t(G)\leq n-2$ if $g(G) = 4$,
    \item$b_t(G)\leq n-2$ if there is a triangle at least one of its vertices is a support vertex in $G$,
    \item $b_t(G)\leq n-1$ if there is a triangle which at least one of its vertices is of degree two in $G$.
\end{enumerate}
    
\end{theorem}

\begin{theorem}\cite{sridharan2007total} If $T$ is a tree on $n$ vertices and $T\not\equiv K_{1,n-1},$ then $$b_t(T)\leq min\{ \Delta(T), \frac{n-1}{3}\}.$$
    
\end{theorem}

In 2014, Rad et al. \cite{rad2014some} proved the following theorem.
\begin{theorem} \cite{rad2014some} For any tree $T$ with maximum degree at least three,
$$b_t(T)\leq \Delta(T)-1.$$
\end{theorem}
In the same paper Rad et al. \cite{rad2014} prove the following two lemmas.
\begin{lemma}
\label{trianglerad}
Let a graph $G$ contains a path $(x,y,z)$ on three vertices such that $\{x,y\}$ are support vertices in $G$, and $y$ is not a support vertex, then $$b_t(G)\leq d(x)+d(y)+d(z)-5.$$
\end{lemma}
Similarly, they prove the following lemma.

\begin{lemma}
\label{4cyclerad}
If a graph $G$ contains a path $(x,y,z)$ on three vertices such that $\{x, y\}$ are not support vertices and $z$ is a support vertex, then
$$b_t(G) \leq d(x) + d(y) + d(z) - 4.$$
\end{lemma}
Also, they prove the following theorem.
\begin{theorem}
\label{4cycle2rad}
Assume that a graph $G$ contains a path $(v_1,v_2,v_3,v_4)$ such that $G-\{v_1,v_2,v_3,v_4\}$ has no isolated vertex. If $G_1$ is the subgraph induced by $v_1,v_2,v_3,v_4$, then $$b_t(G)\leq \sum_{i=1}^{4}d(v_i)-|E(G_1)|-2.$$
\end{theorem}

We use the theorem \ref{borodin} of Borodin at el. \cite{borodin2019all} to prove our theorem. \ref{totalbondagemindegree3planar}

\begin{theorem}\cite{borodin2019all}
\label{borodin}
    Every planar graph with $\delta(G)\geq 3$ has at least one of the following:
    \begin{itemize}
        \item[(a).] a 3-face incident with a $(3,10)$-, $(4,7)$-, or $(5,6)$-edge
        \item[(b).] a 4-face incident either with two $3$-vertices and another $5^-$-vertex or with a $3$-vertex, two $4$-vertices, and a forth vertex of degree at most five;
        \item[(c).] a $5$-face incident with four $3$-vertices and a fifth vertex of degree at most $5$, where all parameters are best possible.
    \end{itemize}
\end{theorem}

Discharging method is a very useful method to prove certain configurations in a planar graph. We use the balanced charging to prove the theorem \ref{girth4totalbondage}. For more details of discharging method, readers can refer the article by Cranston and West \cite{cranston2017introduction}.

\begin{proposition} \cite{cranston2017introduction}
Let $V(G)$ and $F(G)$ be the set of vertices and faces respectively of a planar graph $G$. Denote $l(f)$ to be the length of a face $f$. Then the following qualities hold for $G$.
\begin{equation*}
\begin{aligned}
   &\sum_{\substack{v\in V(G)}}(d(v)-6) + \sum_{\substack{f\in F(G)}}(2l(f)-6) = -12  && \text{vertex charging}\\
   &\sum_{\substack{v\in V(G)}}(2d(v)-6) + \sum_{\substack{f\in F(G)}}(l(f)-6) = -12  && \text{face charging}\\
   &\sum_{\substack{v\in V(G)}}(d(v)-4) + \sum_{\substack{f\in F(G)}}(l(f)-4) = -8  && \text{balanced charging}
\end{aligned}
\end{equation*}
\end{proposition}

\section{Proofs of Lemmas and Theorems}


First, we prove Theorem \ref{sridharanimprovement}.

\begin{proof}
   For this proof, we use the following notations. Let $P=(v,a_1,a_2,\dots, d_{n-1},u)$ be the induced longest path from a support vertex $v$ to a maximum degree vertex $u$ in $T$. Let $T'=T\setminus \{e\}$, where $e$ is an edge in $T$. Define $T_1$ and $T_2$ as the subgraphs of $T'$ that contain $v$ and $u$ respectively. Let $n_1$ and $n_2$ be the respective cardinalities of $V(T_1)$ and $V(T_2)$, and let $E_2\subset E(T_2)$ be a subset of edges such that $|TD(T_2\setminus E_2)|> |TD(T_2)|$. We prove this by induction on the number of vertices of the graph.

First, we consider the base case when $n \leq 7$. Note that $\Delta(T)\geq 3$, so $n\geq 5$. For $n\leq 7$, $\lfloor \frac{n-2}{3}\rfloor=1$. Thus, we want to show that $b_t(T)\leq 1$.

Suppose not, i.e., $b_t(T)\geq 2$, which means deleting one edge does not increase the total domination number. For $n=5$, note that $d(v)=3$, $d(u)=2$, and $d(u,v)=1$. Let $\{e\}=((u,v))$, $2=|TD(T)|=|TD(T\setminus \{e\})|=|TD(T')|=|TD(T_1)|+|TD(T_2)|=4$, a contradiction. So, $b_t(T)\leq 1$ for $n=5$.

Now consider $n=6$. If $d(u,v)=1$, the result follows as $n=5$ by taking $e=\{(u,v)\}$. If $d(u,v)=2$, the result follows similarly as in previous cases by taking $e=\{(v,a_1)\}$.

Now consider $n=7$. If $d(u,v)=1$, the result follows as the previous cases by taking $e=\{u,v\}$. If $d(u,v)=2$, the result follows similarly as in previous cases  by taking $e=\{(v,a_1)\}$. When $d(u,v)=3$, this graph becomes the graph obtained from subdividing a single edge of $K_{1,3}$ three times.

Thus, the result is true for $n\leq 7$. Now assume the result is true for $8\leq k<n$. We want to prove the result is true for $k=n$.

For $n \geq 8$ with $\frac{n-2}{3} \geq 2$, if removing one edge increases the total domination number, the proof is complete. Assume $b_t(T) \geq 2$. For any edge $e \in E(T)$, $|TD(T\setminus \{e\})| = |TD(T)|$. Let $T_1$ and $T_2$ be the resulting subgraphs after deleting the edge $e$ from $T$. If $T_1$ or $T_2$ is isomorphic to $T^1$, a 7-vertex tree obtained by subdividing a single edge of $K_{1,3}$ three times, without loss of generality, suppose $T_1 \equiv T^1$. If $5 \leq |T_2| \leq 7$ and $T_2 \not \equiv T^1$, the result follows from the inductive step. For $|T_2| \leq 3$, it is evident that the removal of at most $2$ edges in $T$ increases the total bondage number. If $|T_2| \geq 7$, then $\frac{n-2}{3} \geq 4$. However, note that in $T^1$, removal of two edges increases the total domination number. Thus, in $T$, removing at most $3$ edges increases the total domination number. Therefore, the result holds true. Assume that for any edge $e \in E(T)$, the resulting sub trees $T_1$ and $T_2$ in $T\setminus \{e\}$ are not isomorphic to $T^1$, a 7-vertex tree obtained by subdividing a single edge of $K_{1,3}$.
three times.

Now, we have to consider three cases. 

\begin{itemize}
    \item[\textbf{Case 1:}] $d(u,v)\geq 3$\\
        Let $e=(a_1,a_2)$ Note that $|TD(T)| =|TD(T\setminus \{e\})|= |TD(T')| = |TD(T_1)| + |TD(T_2)|$. Note that, $|V(T_2)|\leq n-n_1\leq n-3$ (since $n_1\geq 3$). Thus, $b_t(T_2)\leq \frac{n-5}{3}=\frac{n-2}{3}-1$ by the induction hypothesis. So, $|TD(T)| > |TD(T_1)| + |TD(T_2)\setminus E_2| = |TD(T\setminus \{(a_1,a_2)\cup E_2\})|$, and $|\{(a_1,a_2\cup E_2)\}|=\frac{n-2}{3}$. Therefore, the result holds.

    \item[\textbf{Case 2:}] $d(u,v)=2$\\
        \begin{itemize}
            \item[\textbf{Case 2.1:}] $d(u)\geq 4$ or $d(v)\geq 3$\\
                Let $T'= T\setminus \{(v,a_1)\}$ or $T'= T\setminus \{(a_1,u)\}$ for $d(u)\geq 4$ and $d(v)\geq 3$, respectively, so that $|V(T_1)|\geq 3$ and $n_2\leq \frac{n-2}{3}-1$ with $\Delta(T_2)\geq 3$. The result follows by induction, as in Case 1.

            \item[\textbf{Case 2.2:}] $d(u)=3$ and $d(v)=2$\\
                Note that $d(a_1)=2$; otherwise, $b_t(T)\leq 1$ since $|TD(T')| > |TD(T)|$ when $e=(v,a_1)$. Thus, $T$ is a tree obtained from subdividing edges of $K_{1,3}$. Considering $n\geq 8$ and $n\leq 10$, the latter obtained by subdividing each edge of $K_{1,3}$ twice, we can systematically examine cases while preserving $d(u)=3$ and $d(v)=2$. In each case, there exist two edges whose deletion increases the total domination number. Therefore, the result holds.

        \end{itemize}

    \item[\textbf{Case 3:}] $d(u,v)=1$\\
        \begin{itemize}
            \item[\textbf{Case 3.1:}] $d(u)\geq 4$ or $d(v)\geq 3$\\
                Similar to Case 2.1.

            \item[\textbf{Case 3.2:}] $d(u)=3$ and $d(v)=2$\\
                Note that for any $v\in N(u)$ has degree 2; otherwise, the result holds by Case 3.1, and $d(u,v)=1$ for any support vertex $v$. Thus, $|V(T)|\leq 7$, and the result holds by the base case of the induction.

        \end{itemize}
\end{itemize}
\end{proof}




We write a different version of lemma \ref{trianglerad}.
\begin{lemma}
\label{triangle}
Let $G$ be a connected graph containing an induced triangle $(x_1, x_2, x_3)$, where none of the vertices $x_1, x_2, x_3$ are support vertices in $G$. Then, $b_t(G) \leq \sum_{i=1}^{3} d(x_i) - 5$.   
\end{lemma}

Now we prove the lemma \ref{triangle} \begin{proof}.
Let $G \not\equiv C_3$ be a connected graph containing an induced triangle $(x_1,x_2,x_3)$, where none of the vertices $x_1,x_2,x_3$ are support vertices in $G$. We prove by contradiction, assuming $|D'|\leq |D|$ for all $E\subseteq E(G)$ such that $|E|\leq \sum_{i=1}^{3}d(x_i)-5$, where $D$ and $D'$ represent minimal cardinality total dominating sets of $G$ and $G'=G-E$, respectively.

Since $G\not \equiv C_3,$ at least one from $x_1,x_2,x_3$ has degree at least 3. Without loss of generality let $d(x_1)\geq 3$ and  $w\in N(x_1) \setminus \{x_2,x_3\}$. Consider the set $E$ defined by the edges incident with $x_1, x_2, x_3$ except $(x_1,w)$ and $(x_2,x_3)$. Note that $|E|\leq d(x_1) + d(x_2) + d(x_3) - 5$.

Let $D'$ be a total dominating set of $G-E$. Since $x_1$ is a pendant vertex in $G-E$, it follows that $w\in D'$. Additionally, as $(x_2,x_3)$ forms an isolated edge in $G-E$, we have $\{x_2, x_3\}\subset D'$. Now, observe that $D' \setminus \{x_2, x_3\} \cup \{x_1\}$ is a total dominating set of $G$ with size $|D'|- 1< |D'|\leq |D|$, leading to a contradiction.

This contradiction arises regardless of whether $x_1\in D$ or $x_1\not\in D$. Therefore, the result holds.
\end{proof}


We write a different version of lemma for lemma \ref{4cyclerad} and  theorem \ref{4cycle2rad}.
\begin{lemma}
\label{4face}
Let $G$ be a connected graph with minimum degree $\delta(G)\geq 2$, and having an induced 4-cycle $(x_1, x_2, x_3, x_4)$ in $G$. Then,
\[ b_t(G) \leq \sum_{i=1}^{4} d(x_i) - 6. \]
\end{lemma}

Now we prove the lemma \ref{4face} \begin{proof}
Let $G$ be a connected graph with $\delta(G)\geq2$, and having an induced 4-cycle $(x_1, x_2, x_3, x_4)$ in $G$. We prove by contradiction, assuming $|D'|\leq |D|$ for all $E\subseteq E(G)$ such that $|E|\leq \sum_{i=1}^{4}d(x_i)-6$, where $D$ and $D'$ represent minimal cardinality total dominating sets of $G$ and $G'=G-E$, respectively. Let $E$ be the set of edges incident with $x_1, x_2, x_3, x_4$ except $(x_1, x_2)$ and $(x_3, x_4)$. Let $D'$ be a total dominating set of $G'=G-E$. Note that $|E|=\sum_{i=1}^{4}d(x_i)-6$, and $\{x_1, x_2, x_3, x_4\}\subset D'$. Now, $D'\setminus \{x_1, x_4\}$ is a total dominating set of $G$ with size $|D'|-2<|D'| \leq |D|$, a contradiction. Thus, the result holds.
\end{proof}

\begin{lemma}
\label{5face}
Let $G$ be a connected graph with minimum degree $\delta(G)\geq 2$, and having an induced 5-cycle $(x_1, x_2, x_3, x_4, x_5)$ in $G$, where $x_1, x_2, x_3, x_4$ are vertices on the 5-cycle in the clockwise or anticlockwise direction. Then,
\[ b_t(G) \leq \sum_{i=1}^{4} d(x_i) - 5. \]
\end{lemma}
Now we do the proof of lemma \ref{5face} \begin{proof}
Let $G$ be a connected graph with $\delta(G)\geq2$, and having an induced 5-cycle $(x_1, x_2, x_3, x_4, x_5)$ in $G$, where $x_1,x_2,x_3,x_4$ are vertices on the 5-cycle in the clockwise or anticlockwise direction. We prove by contradiction, assuming $|D'|\leq |D|$ for all $E\subseteq E(G)$ such that $|E|\leq \sum_{i=1}^{4}d(x_i)-5$, where $D$ and $D'$ represent minimal cardinality total dominating sets of $G$ and $G'=G-E$, respectively and $\{x_1, x_2, x_3, x_4\}$ is a set of vertices on the 5-face in the clockwise or anticlockwise direction. 
\vspace{0.3cm}
Now consider two cases.
\begin{itemize}
    \item[(a).] $d(x_5)\geq 3$
    \item[(b).] $d(x_5)=2$
\end{itemize}
In the case (a), let $E$ be the set of edges incident with $x_1, x_2, x_3, x_4$ except $(x_1, x_2)$ and $(x_3, x_4)$. Let $D'$ be a total dominating set of $G'=G-E$. Note that $|E|=\sum_{i=1}^{4}d(x_i)-5$ and $\{x_1, x_2, x_3, x_4\}\subset D'$. Now, $D'\setminus \{x_1,x_4\}$ is a total dominating set of $G$ with $|D'|-2<|D'|\leq |D|$, a contradiction.

In case (b), let $E$ be the set of edges incident with $x_1, x_2, x_3, x_4$ except $(x_1,x_2), (x_3,x_4),$ and $(x_4,x_5)$. Note that $x_4$ is a pendant vertex in $G'=G-E$. $\{x_1, x_2, x_4\}\subset D'$. Now consider two sub-cases.
\begin{itemize}
    \item $x_5\in D'$
    \item $x_5 \not \in D'$
\end{itemize}
In the former case, $D'\setminus \{x_2\}$ is a total dominating set of $G$ with $|D'|-1<|D'|\leq |D|$, a contradiction. In the latter case, $x_3\in D'$, and thus $D'\setminus \{x_1\}$ is a total dominating set of $G$ with $|D'|-1<|D'|\leq |D|$, a contradiction.

Thus, the result holds.
\end{proof}


Now we prove the theorem \ref{totalbondagemindegree3planar} \begin{proof}.
Let $G$ be a connected planar graph with $\delta(G)\geq 3$. Suppose, for the sake of contradiction that $b_t(G)\geq \Delta(G) + 9$ and $b_t(G)\geq 11$.

By Theorem \ref{borodin}, $G$ has one of the following cases:

\begin{itemize}
\item[(a)] a 3-face incident with a $(3,10)$-, $(4,7)$-, or $(5,6)$-edge
\item[(b)] a 4-face incident either with two $3$-vertices and another $5^-$-vertex or with a $3$-vertex, two $4$-vertices, and a fourth vertex of degree at most five
\item[(c)] a 5-face incident with four $3$-vertices and a fifth vertex of degree at most $5$, where all parameters are best possible.
\end{itemize}

\textbf{Case (a):}
If $G$ has a 3-face incident with a $(3,10)$-, $(4,7)$-, or $(5,6)$-edge, then by Theorem \ref{triangle}, $G$ has a triangle $(x_1,x_2,x_3)$ with a degree sum at most $13+ \Delta$. Applying Theorem \ref{triangle}, we get $b_t(G)\leq \Delta + 8$, which leads to a contradiction.

\textbf{Case (b):}
If $G$ has a 4-face incident as described in case (b), then $G$ has a 4-face $(x_1,x_2,x_3,x_4)$ with a degree sum at most $\Delta +11$ or $16$. Applying Theorem \ref{4face}, we find that $b_t(G)\leq \Delta +5$, leading to a contradiction.

\textbf{Case (c):}
If $G$ has case (c), then $G$ has a 5-face with the degree sum of 4-vertices incident with the 5-face is at most $12$. By Theorem \ref{5face}, $b_t(G)\leq 7$, which is a contradiction.

Therefore, in all cases, we reach a contradiction. Thus, the result holds.
\end{proof}



Proof of theorem \ref{configurationgirth4}. \begin{proof}
Let $G$ be a connected planar graph with $\delta(G) \geq 3$ and girth $g(G)\geq 4$ without any of the given configurations. We will proceed by balanced charging. Every face $f$ begins with charge $l(f)-4$, and every vertex begins with charge $d(v)-4$. Hence, only $3^-$-vertices start with an initial negative charge. All $4^+$-vertices have a non-negative charge. We begin with the following discharging rule:
\begin{itemize}
    \item[(R)] Every $3$-vertex takes $\frac{1}{3}$ from each adjacent vertex.
\end{itemize}
We want to show that every vertex and face ends with a non-negative charge. Suppose not. We begin by checking the vertices. Let $v$ be a $3$-vertex. Note that $v$ starts with charge $-1$. Thus, $v$ receives less than $1$. By $(R)$, $v$ receives a total of $3*\frac{1}{3}=1$, a contradiction. Now consider a 4-vertex $v$. Note that we avoid $(3,4^-)$-edges in $G$, thus 4-vertices do not lose charge to 3-vertices, so 4-vertices start with charge $0$ and remain non-negative, a contradiction. Now consider a $5$-vertex $v$. $v$ starts with charge $1$ and thus loses more than $1$. $v$ loses charge to incident $3$-vertices. If $|N(v)\cap S_3|\leq 3$, then $v$ loses at most $3*\frac{1}{3}=1$, and stays non-negative, a contradiction. Thus, $|N(v)\cap S_3|\geq 4$. But we avoid the configuration $(b)$ in $G$, thus $v$ ends with a non-negative charge, a contradiction. Now consider a $6^+$-vertex $v$. Note that $v$ starts with charge $d(v)-4$, and thus $v$ loses more than $d(v)-4$. $v$ loses the highest charge when $|N(v)\cap S_3|=d(v)$. So, assume $|N(v)\cap S_3|=d(v)$. Now, $v$ loses $d(v)*\frac{1}{3}$ and retains $d(v)-4-d(v)*\frac{1}{3}=\frac{2}{3}d(v)-4<0$. This implies that $d(v)<6$, a contradiction. Thus, every vertex ends with a non-negative charge.

Now consider the faces. Note that every face $f$ starts with charge $l(f)-4$ and does not lose charge with the discharging rule. Suppose $f$ ends with negative charge. Thus, $l(f)-4<0$, which implies that $l(f)<4$, a contradiction since $g(G)\geq 4$. Thus, every face ends with a non-negative charge.

This leads to a contradiction, as we started with a graph with a non-negative charge. Thus, the result holds.
\end{proof}

Now, we will prove the lemma \ref{3vertices at distance 2}.\\ Let $G$ be a graph with $\delta(G)\geq 3$, and $\{u_1,u_2\}\subset V(G)$ such that $d(u_i)=3$ for $i\in \{1,2\}$, and $d(u_1,u_2)=2$. Suppose to the contrary that $|D'|\leq |D|$ for all $E\subset E(G)$ such that $|E|\leq \Delta +3$, where $D$ and $D'$ represent minimal cardinality total dominating sets of $G$ and $G'=G-E$, respectively. Let $v\in N(u_1)\cap N(u_2)$, and $E$ be the set of edges incidents with $u_1,u_2,u_2'$ except $(u_1,v), (u_2,u_2')$, where $u_2'\in N(u_2)\setminus \{v\}$. Note that $|E|\leq \Delta +3$. Now, $v\in D'$ as $u_1$ is a pendant vertex in $G-E$. Also, $\{u_2,u_2'\}\subset D'$. Note that, $D'\setminus\{u_2'\}$ is a total dominating set of $G$ of size $|D'|-1<|D'|\leq |D|$, a contradiction. Thus, the result holds.

\vspace{0.3cm}
Now we prove the theorem \ref{girth4totalbondage}.
\begin{proof}
    Suppose $G$ be a connected planar graph $G$ with $\delta(G) \geq 3$ and girth $g(G)\geq 4$. Then by theorem \ref{configurationgirth4}, $G$ has 
    \begin{enumerate}
   \item[(a)] $(3,4^-)$-edge or
   \item[(b)] $v \in S_5$ such that $|N(v)\cap S_3|\geq 4$.
\end{enumerate}
In the theorem \ref{girth4totalbondage}, we avoid the case $(a)$, where $G$ having an $(3,4^-)$-edge. Thus, we only need to consider the case that $G$ has the configuration $(b)$. Now, we have two 3-vertices at a distance of exactly 2. Thus by lemma \ref{3vertices at distance 2}, the result holds.  
\end{proof}

\vspace{0.3cm}

Now we prove the theorem \ref{distant2totalbondage}.
\begin{proof}
 
  Let $G$ be a connected graph with $\delta(G)\geq 2$, and $\{u_1,u_2\}\subset V(G)$ such that $d(u_i)=2$ for $i\in \{1,2\}$, and $d(u_1,u_2)\leq 3$. Suppose to the contrary that $|D'|\leq |D|$ for all $E\subset E(G)$ such that $|E|\leq \Delta +1$, where $D$ and $D'$ represent minimal cardinality total dominating sets of $G$ and $G'=G-E$, respectively. 
  We have three cases:
  \begin{enumerate}
      \item[(a)] $d(u_1,u_2)=3$
      \item[(b)]$d(u_1,u_2)=2$
      \item[(c)]$d(u_1,u_2)=1$
  \end{enumerate}
  
  \textbf{Case (a):} Let $u_1vwu_2$ be a path containing $u_1,u_2$ with $d(u_1,u_2)=3$, where $\{v,w\}\subset V(G)$. Now, let $E$ be the set of edges incidents with $\{u_1,v,u_2\}$ except $(u_1,v)$ and $(w,u_2)$. Note that, $|E|\leq \Delta +1$, and $\{u_1,v,w\}\subset D'$. Now, $D'\setminus \{u_1\}$ is a total dominating set of $G$ of size $|D'|-1<|D'|\leq |D|$, a contradiction.\\
  \textbf{Case (b):} Let $vu_1wu_2$ be a path containing $u_1,u_2$ with $d(u_1,u_2)=2$, where $\{v,w\}\subset V(G)$. Now, let $E$ be the set of edges incidents with $\{v,u_1,u_2\}$ except $(v,u_1)$ and $(w,u_2)$. Note that, $|E|\leq \Delta+1$, and $\{v,u_1,w\}\subset D'$. Now, $D'\setminus \{v\}$ is a total dominating set of $G$ of size $|D'|-1<|D'|\leq |D|$, a contradiction.\\
  \textbf{Case (c):} Let $vu_1u_2$ be a path containing $u_1,u_2$ with $d(u_1,u_2)=1$, where $\{v\}\in V(G)$. Now, let $E$ be the set of edges incidents with $\{v',u_1,u_2\}$ except $(u_1,u_2)$  and $(v,v')$, where $v'\in N(v)\setminus \{u_1\}$. Note that $|E|\leq \Delta(G)+1$, and $\{v,u_1,u_2\}\subset D'$. Now, $D'\setminus \{u_2\}$ is a total dominating set of $G$ of size $|D'|-1<|D'|\leq |D|$, a contradiction. 
  Thus, the result holds.

\end{proof}

Now, we will prove theorem \ref{completemultipartitegraphtotalbondate}.

\begin{proof}
Let $K_{n_1,n_2,\dots, n_k}$ be a complete multipartite graph with $2\leq n_k\leq n_{k-1}\leq \dots \leq n_2\leq n_1 $.
Suppose $k=2$, then $K_{n_1,n_2}=n_2$ by Theorem \ref{kullibipartite}. Let $A_i$ be the vertex sets of $K_{n_1,n_2,\dots, n_k}$ with $|A_i|=n_i$ for $i \in \{1,2,\dots, k\}$. Choose distinct vertices $u^1_1,u^1_2$ from $A_1$, and $u^k_1, u^k_2$ from $A_k$. 
Let $E$ be the set of all edges incident with $u^1_1$ and $u^1_2$ except $(u^1_1,u^k_1)$ and $(u^1_2,u^k_2)$. Observe that $|E|=2n-2n_1-2$, and $u^k_1,u^k_2$ are pendant vertices in $K_{n_1,n_2,\dots, n_k}\setminus E$. 
Thus, $u^k_1,u^k_2$ are in a minimum cardinality total dominating set $D$ of $K_{n_1,n_2,\dots, n_k}\setminus E$. Since there is no edge between $u^k_1$ and $u^k_2$, $u_i \in D$ for $u_i \in (A_1\cup A_2\cup \dots \cup A_k\setminus \{u^1_1,u^1_2,u^k_1,u^k_2\})$. Therefore, $2=|TD(K_{n_1,n_2,\dots, n_k})|>|D|=|TD(K_{n_1,n_2,\dots, n_k}\setminus E)|$. Thus, the result holds.

\end{proof}



\end{document}